\input amstex
\documentstyle{amsppt}
%
%
\nopagenumbers
\def\const{\operatorname{const}}
\def\negskp{\hskip -2pt}
\pagewidth{360pt}
\pageheight{606pt}
\rightheadtext{Dynamical systems admitting normal blow-up of points.}
\topmatter
\title
NEWTONIAN DYNAMICAL SYSTEMS ADMITTING
NORMAL BLOW-UP OF POINTS.
\endtitle
\author
Ruslan~A\.~Sharipov
\endauthor
\abstract
Class of Newtonian dynamical systems admitting normal blow-up of
points in Riemannian manifolds is considered. Geometric interpretation
for weak normality condition, which arose earlier in the theory of
dynamical systems admitting the normal shift of hypersurfaces, is found.
\endabstract
\address Rabochaya~str\.~5, 450003, Ufa, Russia
\endaddress
\email \vtop to 20pt{\hsize=280pt\noindent
R\_\hskip 1pt Sharipov\@ic.bashedu.ru\newline
ruslan-sharipov\@usa.net\vss}
\endemail
\urladdr
http:/\negskp/www.geocities.com/CapeCanaveral/Lab/5341
\endurladdr
\endtopmatter
\loadbold
\document
\head
1. Introduction.
\endhead
     Let $S$ be a hypersurface in Riemannian manifold $M$. One of the
ways for deforming $S$ consists in shifting points, which constitute $S$,
along trajectories of some Newtonian dynamical system. Such situation
arises in describing the propagation of electromagnetic wave (light) in
non-homogeneous media in the limit of geometric optics. Hypersurface
$S$ models wave front set (the set of points with constant phase),
while trajectories of shift model light beams. Newtonian dynamics of
points of Riemannian manifold $M$ in local coordinates $x^1,\,\ldots,
\,x^n$ in $M$ is described by a system of $n$ ordinary differential
equations
$$
\ddot x^k+\sum^n_{i=1}\sum^n_{j=1}\Gamma^k_{ij}\,\dot x^i\,\dot x^j
=F^k(x^1,\ldots,x^n,\dot x^1,\ldots,\dot x^n),\hskip -2em
\tag1.1
$$
where $k=1,\,\ldots,\,n$. Here $\Gamma^k_{ij}$ are components of metric
connection for basic metric $\bold g$ of the manifold $M$. Quantities
$F^k$ are components of force vector $\bold F$. They determine force
field of dynamical system \thetag{1.1}. In the equations \thetag{1.1}
they play role of perturbing factor, due to them trajectories of
dynamical system \thetag{1.1} differ from that of geodesic flow for
the metric $\bold g$.\par
     At each point $p$ of hypersurface $S$ we fix some vector of
initial velocity $\bold v(p)$ and determine trajectories coming
out from all points of hypersurface $S$ by setting the following
Cauchy problem for the equations \thetag{1.1}:
$$
\xalignat 2
&\quad x^k\,\hbox{\vrule height 8pt depth 8pt width 0.5pt}_{\,t=0}
=x^k(p),
&&\dot x^k\,\hbox{\vrule height 8pt depth 8pt width 0.5pt}_{\,t=0}=
v^k(p).\hskip -2em
\tag1.2
\endxalignat
$$
Here $v^k(p)$ are the components of vector $\bold v(p)$ in local
coordinates $x^1,\,\ldots,\,x^n$. Having displaced for time $t$
along trajectories determined by initial data \thetag{1.2}, points
of the hypersurface $S$ constitute another hypersurface $S_t$. As
a result we obtain a family of hypersurfaces and a family of shift
maps $f_t\!:S\to S_t$, which are local diffeomorphisms for
sufficiently small values of parameter $t$. All this family of maps
is called {\bf a construction of shift} or simply {\bf a shift} of
hypersurface $S$ along trajectories of dynamical system \thetag{1.1}.
\definition{Definition 1.1} Shift $f_t\!:S\to S_t$ of hypersurface
$S$ along trajectories of Newtonian dynamical system with force
field $\bold F$ is called {\bf a normal shift} if hypersurfaces
$S_t$ are orthogonal to the trajectories of shift.
\enddefinition
     In order to construct the normal shift we should, at least,
choose initial velocities $\bold v(p)$ being perpendicular to
initial hypersurface $S$, i\.~e\. $\bold v(p)=\nu(p)\cdot\bold n(p)$,
where $\bold n(p)$ is a unitary normal vector to $S$ at the point
$p$:
$$
\xalignat 2
&\quad x^k\,\hbox{\vrule height 8pt depth 8pt width 0.5pt}_{\,t=0}
=x^k(p),
&&\dot x^k\,\hbox{\vrule height 8pt depth 8pt width 0.5pt}_{\,t=0}=
\nu(p)\cdot n^k(p).\hskip -2em
\tag1.3
\endxalignat
$$
But this is not sufficient. Initial data \thetag{1.3} by themselves
do not provide orthogonality of $S_t$ and trajectories of shift for
$t\neq 0$. We are to use other opportunities due to the choice of
hypersurface $S$, choice of modulus of initial velocity $\nu(p)$ on
$S$, and choice of the force field $\bold F$ of dynamical system
\thetag{1.1}. If we choose hypersurface $S$, then, in the case
of success, we would have the construction of normal shift realized
only on some special hypersurface $S$ (or in some special class of
hypersurfaces). In paper \cite{1} we left the choice of $S$ to be
arbitrary, and have concentrated efforts to the choice of function
$\nu(p)$ in \thetag{1.3}. Then we found that the proper choice of
$\nu(p)$ on an arbitrary hypersurface $S$ is possible only under
some definite restrictions for the choice of force field of dynamical
system \thetag{1.1}. This became the origin for the theory of dynamical
systems {\bf admitting the normal shift}. It was developed in the series
of papers \cite{1--16}. On the base of these papers two theses were
prepared: thesis for the degree of {\bf Doctor of Sciences in Russia}
\cite{17} and thesis for the degree of {\bf Candidate of Sciences in
Russia} \cite{18}.\par
    In \cite{1--18} we restricted ourselves to the case of smooth
hypersurfaces $S_t$ with no singular points. However, in the process
of shifting $f_t\!:S\to S_t$ some singular points can appear (they
are called {\bf caustics}). In particular, we can observe the collapse
of $S_t$ into a point at some instant of time $t=t_0$ followed by a
blow-up of this point into further series of smooth hypersurfaces for
$t>t_0$. Without loss of generality we can assume that $t_0=0$.
Then initial hypersurface $S=\{p_0\}$ consisting of only one point
$p_0$ appears to be singular, and for $t>0$ we have blow-up of this
point $p_0$ into a series of smooth hypersurfaces $S_t$.
\definition{Definition 1.2} Blow-up $f_t\!:S\to S_t$ of singular
one-point hypersurface $S=\{p_0\}$ along trajectories of Newtonian
dynamical system with force field $\bold F$ is called {\bf a normal
blow-up} if for $t>0$ smooth hypersurfaces $S_t$ are orthogonal
to the trajectories of this blow-up.
\enddefinition
     The idea to consider blow-ups of one-point sets in the framework
of normal shift was suggested by A.~V.~Bolsinov and A.~T.~Fomenko
when author was reporting results of thesis \cite{17} in a seminar
at Moscow State University. Partially this idea was realized in
\cite{19}. In that paper was shown that Newtonian dynamical systems
admitting the normal shift of hypersurfaces are able to implement
{\bf normal blow-up} of any point $p_0$ in Riemannian manifold $M$.
More completely the idea of A.~V.~Bolsinov and A.~T.~Fomenko can
be realized in special investigation. This is the main goal of
present paper.
\head
2. Geometry of normal blow-up.
\endhead
\parshape 18 0pt 360pt 0pt 360pt 0pt 360pt 160pt 200pt 160pt 200pt
160pt 200pt 160pt 200pt 160pt 200pt 160pt 200pt 160pt 200pt
160pt 200pt 160pt 200pt 160pt 200pt 160pt 200pt 160pt 200pt
160pt 200pt 160pt 200pt 0pt 360pt
     Let $p_0$ be some point of Riemannian manifold. Let's consider
normal blow-up of this point $f_t\!:\{p_0\}\to S_t$ along trajectories
of Newtonian dynamical system \thetag{1.1}. Consider some particular
trajectory of shift. For $t=0$ it passes through the point $p_0$.
\vadjust{\vskip 10pt\hbox to 0pt{\kern 5pt\hbox{\special{em:graph
pst-02a.gif}}\hss}\vskip -10pt}Let $\bold v(0)$ be the velocity
vector corresponding to the time instant $t=0$. If $\bold v(0)\neq 0$,
then this vector can be normalized to unit length:
$$
\bold n=\frac{\bold v(0)}{|\bold v(0)|}.\hskip -2em
\tag2.1
$$
Unit vectors \thetag{2.1} for various trajectories belong to the same
tangent space $T_{p_0}(M)$. They can be identified with radius-vectors
of points on unit sphere $\sigma$ in the space $T_{p_0}(M)$.
\definition{Definition 2.1}\parshape 3 160pt 200pt 160pt 200pt
0pt 360pt                             
The blow-up $f_t\!:\{p_0\}\to S_t$ of
the $p_0$ along trajectories of Newtonian dynamical system \thetag{1.1}
is called {\bf regular} if velocity vectors $\bold v(0)$ at the point
$p_0$ are non-zero for all trajectories of this blow-up and if
points corresponding to unit vectors \thetag{2.1} fill the whole
surface of unit sphere $\sigma$ in $T_{p_0}(M)$.
\enddefinition
     In the case of regular blow-up all hypersurfaces $S_t$ possess
spherical topology for sufficiently small values of parameter $t\neq 0$.
Points $q$ of the unit sphere $\sigma$ in $T_{p_0}(M)$ can be used to
parameterize points of hypersurface $S_t$. In order to do it we shall
write initial data determining regular blow-up of the point $p_0$ as
follows:
$$
\xalignat 2
&\quad x^k\,\hbox{\vrule height 8pt depth 8pt width 0.5pt}_{\,t=0}
=x^k(p_0),
&&\dot x^k\,\hbox{\vrule height 8pt depth 8pt width 0.5pt}_{\,t=0}=
\nu(q)\cdot n^k(q)\hskip -2em
\tag2.2
\endxalignat
$$
Trajectories of Newtonian dynamical system \thetag{1.1} fixed by initial
data \thetag{2.2} determine a family of maps $f_t\!:\sigma\to S_t$ being
diffeomorphisms for sufficiently small values of parameter $t\neq 0$.
\par
    Suppose that we have regular blow-up of the point $p_0$ being normal
in the sense of definition~1.2. Let's consider the solution of Cauchy
problem \thetag{2.2} for the equations \thetag{1.1}. This is the set of
$n$ functions $x^1(q,t),\,\ldots,\,x^n(q,t)$. Due to initial data
\thetag{2.2} we can write Taylor expansions for these functions at the
point $t=0$:
$$
x^k(q,t)=x^k(p_0)+\nu(q)\,n^k(q)\cdot t+O(t).\hskip -2em
\tag2.3
$$
Denote $\bold v(q)=\nu(q)\cdot\bold n(q)$. Vector $\bold v(q)\in T_{p_0}(M)$
has the meaning of initial velocity for trajectory that corresponds to the
point $q$ on unit sphere $\sigma$. In terms of components of vector $\bold
v(q)$ the expansions \thetag{2.3} can be rewritten as
$$
x^k(q,t)=x^k(p_0)+v^k(q)\cdot t+O(t).\hskip -2em
\tag2.4
$$
Let $u^1,\,\ldots,\,u^{n-1}$ be local coordinates of the point $q$ on
unit sphere $\sigma$ in $T_{p_0}(M)$. Due to local diffeomorphisms of
blow-up $f_t\!:\sigma\to S_t$ they can be used as local coordinates on
hypersurfaces $S_t$. Let's represent functions $x^k(q,t)$ and their
expansions \thetag{2.4} in local coordinates $u^1,\,\ldots,\,u^{n-1}$:
$$
x^k(u^1,\ldots,u^{n-1},t)=x^k(p_0)+v^k(u^1,\ldots,u^{n-1})\cdot t+O(t).
\hskip -2em
\tag2.5
$$
Time derivatives of the functions \thetag{2.5} determine velocity vector
on trajectories of blow-up, their derivatives in $u^1,\,\ldots,\,u^{n-1}$
determine tangent vectors to $S_t$:
$$
\xalignat 2
&\bold v=\sum^n_{k=1}\frac{\partial x^k}{\partial t}\cdot\frac{\partial}
{\partial x^k},
&&\boldsymbol\tau_i=\sum^n_{k=1}\frac{\partial x^k}{\partial u^i}\cdot
\frac{\partial}{\partial x^k}.\hskip -2em
\tag2.6
\endxalignat
$$
Change in any one of parameters $u^1,\,\ldots,\,u^{n-1}$ leads to
transfer from one trajectory of blow-up to another. Therefore
vectors $\boldsymbol\tau_1,\,\ldots,\,\boldsymbol\tau_{n-1}$ are
called {\bf vectors of variation of trajectories} or simply
{\bf vectors of variation}. For the components of the vector
$\bold v(t)$ and for the components of vectors $\boldsymbol
\tau_i(t)$ from \thetag{2.5} we derive 
$$
\align
&v^k(u^1,\ldots,u^{n-1},t)=\frac{\partial x^k}{\partial t}=
v^k(u^1,\ldots,u^{n-1})+O(1),\hskip -2em
\tag2.7\\
\vspace{1ex}
&\tau^k_i(u^1,\ldots,u^{n-1},t)=\frac{\partial x^k}{\partial u^i}=
\frac{\partial v^k(u^1,\ldots,u^{n-1})}{\partial u^i}\cdot t+O(t).
\hskip -2em
\tag2.8
\endalign
$$
The normality condition for the blow-up $f_t\!:\{p_0\}\to S_t$
implies orthogonality of velocity vector $\bold v$ to all vectors
$\boldsymbol\tau_1,\,\ldots,\,\boldsymbol\tau_{n-1}$ in \thetag{2.6}.
Let's write this condition as the condition of vanishing of scalar
products $(\bold v\,|\,\boldsymbol\tau_i)$:
$$
(\bold v\,|\,\boldsymbol\tau_i)=\sum^n_{k=1}\sum^n_{r=1}
g_{kr}(x^1,\ldots,x^n)\,\frac{\partial x^k}{\partial t}
\frac{\partial x^r}{\partial u^i}=0.\hskip -2em
\tag2.9
$$
Substituting the expansions \thetag{2.5}, \thetag{2.7}, and
\thetag{2.8} into the equality \thetag{2.9}, we determine the
asymptotics of left hand side of this equality as $t\to 0$:
$$
\sum^n_{k=1}\sum^n_{r=1}g_{kr}(p_0)\,v^k(u^1,\ldots,u^{n-1})\,
\frac{\partial v^k(u^1,\ldots,u^{n-1})}{\partial u^i}\cdot t+
O(t)=0.\hskip -2em
\tag2.10
$$
Right hand side of \thetag{2.10} is identically zero. Therefore
from \thetag{2.10} we get
$$
\sum^n_{k=1}\sum^n_{r=1}g_{kr}(p_0)\,v^k(u^1,\ldots,u^{n-1})\,
\frac{\partial v^k(u^1,\ldots,u^{n-1})}{\partial u^i}=0.
\hskip -2em
\tag2.11
$$
Here $g_{kr}(p_0)=g_{kr}(x^1(p_0),\ldots,x^n(p_0))$ are the components
of metric tensor at the point $p_0$ referred to local coordinates
$x^1,\,\ldots,\,x^n$ in $M$. Looking attentively at the left hand side
of \thetag{2.11}, we see that it is exactly the scalar product of the
vector of initial velocity $\bold v(q)$ and the derivative of this vector
with respect to parameter $u^i$:
$$
(\bold v(q)\,|\,\bold v^{\,\prime}_{u^i}(q))=
\frac{1}{2}\,\frac{\partial |\bold v(q)|^2}{\partial u^i}=
\nu(q)\,\frac{\partial\nu(q)}{\partial u^i}=0.
\hskip -2em
\tag2.12
$$
This equality \thetag{2.12} proves the following theorem.
\proclaim{Theorem 2.1} For the regular blow-up $f_t\!:\{p_0\}\to S_t$
of the point $p_0$ along the trajectories of Newtonian dynamical
system \thetag{1.1} to be normal it should be determined by initial
data \thetag{2.2} with $\nu(q)$ being constant: \pagebreak
$\nu(q)=\const\neq 0$.
\endproclaim
     Let's denote by $\nu_0$ the constant appeared in theorem~2.1.
Then we can write the initial data \thetag{2.2} in the following form:
$$
\xalignat 2
&\quad x^k\,\hbox{\vrule height 8pt depth 8pt width 0.5pt}_{\,t=0}
=x^k(p_0),
&&\dot x^k\,\hbox{\vrule height 8pt depth 8pt width 0.5pt}_{\,t=0}=
\nu_0\cdot n^k(q)\hskip -2em
\tag2.13
\endxalignat
$$
Further we state the following definition, which is central for the
theory developed in present paper below.
\definition{Definition 2.2} Newtonian dynamical system \thetag{1.1}
with force field $\bold F$ on Riemannian manifold $M$ is called a system
{\bf admitting normal blow-up of points} if for any point $p_0\in M$
and for arbitrary positive constant $\nu_0>0$ initial data \thetag{2.13}
determine the normal blow-up $f_t\!:p_0\to S_t$ of the point $p_0$ along
trajectories of this dynamical system.
\enddefinition
    Definition~2.2 was first formulated in paper \cite{19}. Theorem~2.1
shows that the condition $\nu(q)=\nu_0=\const$ built into definition~2.1
is absolutely inevitable.\par
    Let's substitute $\nu(q)=\nu_0$ into the expansions \thetag{2.5},
\thetag{2.7}, and \thetag{2.8} and take into account the above notation
$\bold v(q)=\nu(q)\cdot\bold n(q)$. As a result we get the expansions
$$
\align
&x^k(u^1,\ldots,u^{n-1},t)=x^k(p_0)+\nu_0\,n^k(u^1,\ldots,u^{n-1})
\cdot t+O(t),\hskip -2em
\tag2.14\\
\vspace{2ex}
&v^k(u^1,\ldots,u^{n-1},t)=\nu_0\,n^k(u^1,\ldots,u^{n-1})+O(1),
\hskip -2em
\tag2.15\\
\vspace{1ex}
&\tau^k_i(u^1,\ldots,u^{n-1},t)=\nu_0\,\frac{\partial n^k(u^1,\ldots,
u^{n-1})}{\partial u^i}\cdot t+O(t),\hskip -2em
\tag2.16
\endalign
$$
which hold due to initial data \thetag{2.13} determining blow-up of the
point $p_0$ along trajectories of dynamical system \thetag{1.1}.
\head
3. Dynamical systems\\
admitting normal blow-up of points.
\endhead
    Definition~2.2 introduces new special class of Newtonian dynamical
systems. According to the results of \cite{19}, it is not empty (see
theorem~12.1 in \cite{19}). In present paper we study this new class
of dynamical systems introduced by definition~2.2.\par
    Let $\bold F$ be the force field of Newtonian dynamical system
admitting normal blow-up of points. Then, according to definition~2.2,
by choosing an arbitrary point $p_0\in M$ and by fixing an arbitrary
constant $\nu_0>0$ one can construct normal blow-up $f_t\!:\{p_0\}\to
S_t$. In order to study this blow-up we consider hypersurfaces $S_t$
with spherical topology, determine local coordinates $u^1,\,\ldots,
\,u^{n-1}$ transferred from unit sphere $\sigma$ to $S_t$, and
define vectors \thetag{2.6} on trajectories of this blow-up. Then we
introduce the following scalar products:
$$
\varphi_i=(\bold v\,|\,\boldsymbol\tau_i).\hskip -2em
\tag3.1
$$
Such scalar products were already considered above in formula
\thetag{2.9}. In thesis \cite{17} they were called {\bf the functions
of deviation}. Functions of deviations \thetag{3.1} are the measure of
deviation of blow-up $f_t\!:\{p_0\}\to S_t$ from normality. In the
case of normal blow-up all these functions are identically zero:
$\varphi_i=0$.\par
    Vanishing of the functions of deviation $\varphi_i$ at the initial
instant of time $t=0$ follows from initial conditions \thetag{2.13}
regardless to the choice of force field $\bold F$ of Newtonian dynamical
system \thetag{1.1}:
$$
\varphi_i\,\hbox{\vrule height 8pt depth 8pt width 0.5pt}_{\,t=0}
=\lim_{t\to 0}\varphi_i(u^1,\ldots,u^{n-1},t)=0.\hskip -2em
\tag3.2
$$
Indeed, as $t\to 0$ vector of velocity tends to its limit value $\bold
v(0)=\nu_0\cdot\bold n(q)$, while vector $\boldsymbol\tau_i$ tends to
zero, this follows from the expansions \thetag{2.8} for its components.
Hence scalar product $\varphi_i=(\bold v\,|\,\boldsymbol\tau_i)$ tends
to zero.\par
    Apart from \thetag{3.2}, identical vanishing of the functions of
deviation in the case of normal blow-up implies vanishing of their
time derivatives $\dot\varphi_i$:
$$
\dot\varphi_i\,\hbox{\vrule height 8pt depth 8pt width 0.5pt}_{\,t=0}
=\lim_{t\to 0}\frac{\partial\varphi_i(u^1,\ldots,u^{n-1},t)}
{\partial t}=0.\hskip -2em
\tag3.3
$$
In calculating $\dot\varphi_i$ for $t\neq 0$ we can replace differentiation
in $t$ by the covariant differentiation with respect to parameter $t$ along
the trajectory of blow-up:
$$
\dot\varphi_i=\nabla_t\varphi_i=\nabla_t(\bold v\,|\,\boldsymbol
\tau_i)=(\nabla_t\bold v\,|\,\boldsymbol\tau_i)+(\bold v\,|\,
\nabla_t\boldsymbol\tau_i).\hskip -2em
\tag3.4
$$
For $\nabla_t\bold v$ we have $\nabla_t\bold v=\bold F$. This follows
from the equations of Newtonian dynamics \thetag{1.1}. Therefore formula
\thetag{3.4} for $\dot\varphi_i$ now is written as follows:
$$
\dot\varphi_i=(\bold F\,|\,\boldsymbol\tau_i)+(\bold v\,|\,
\nabla_t\boldsymbol\tau_i).\hskip -2em
\tag3.5
$$
Vector $\bold F$ has a finite limit as $t\to 0$, it is determined by
initial conditions \thetag{2.13}: $\bold F\to\bold F(p_0,\nu_0\cdot
\bold n(q))$. While vector $\tau_i$ tends to zero. Therefore first
summand in \thetag{3.5} vanishes in the limit as $t\to 0$. Now consider
vector $\nabla_t\boldsymbol\tau_i$ in the second summand. Let's write
components of this vector:
$$
\nabla_t\tau^k_i=\frac{\partial\tau^k_i}{\partial t}+\sum^n_{r=1}
\sum^n_{s=1}\Gamma^k_{rs}\,v^r\,\tau^s_i.\hskip -2em
\tag3.6
$$
Then let's calculate limit as $t\to 0$ in formula \thetag{3.6},
using the expansions \thetag{2.14}, \thetag{2.15}, and \thetag{2.16}
for this purpose. As a result we get
$$
\lim_{t\to 0}\nabla_t\tau^k_i=\nu_0\,\frac{\partial
n^k(u^1,\ldots,u^{n-1})}{\partial u^i}.\hskip -2em
\tag3.7
$$
Here $n^k(u^1,\ldots,u^{n-1})$ are components of unitary vector
$\bold n(q)$ being the radius vector of the point $q$ on unit sphere
$\sigma$ in tangent space $T_{p_0}(M)$. Derivatives of the vector
$\bold n(u^1,\ldots,u^{n-1})$ in $u^1,\,\ldots,\,u^{n-1}$ are coordinate
tangent vectors to the sphere $\sigma$ in local coordinates $u^1,\,\ldots,
\,u^{n-1}$. Let's denote these vectors by $\bold K_1,\,\ldots,\,\bold
K_{n-1}$:
$$
\bold K_i(q)=\sum^n_{k=1}\frac{\partial\bold n^k(u^1,\ldots,u^{n-1})}
{\partial u^i}\cdot\frac{\partial}{\partial x^k}.\hskip -2em
\tag3.8
$$
Then the relationship \thetag{3.7} can be rewritten as
$$
\nabla_t\boldsymbol\tau_i\to\nu_0\cdot\bold K_i(q)\text{\ \ as \ }
t\to 0.\hskip -2em
\tag3.9
$$
For vector of velocity $\bold v$, as was mentioned above, we
have the relationship
$$
\bold v\to\nu_0\cdot\bold n(q)\text{\ \ as \ }t\to 0.\hskip -2em
\tag3.10
$$
From \thetag{3.9} and \thetag{3.10} we derive vanishing of the second
summand in right hand side of \thetag{3.5} as $t\to 0$. Indeed, vector
$\bold K_i(q)$ is tangent to unit sphere $\sigma$ at the point $q$,
while vector $\bold n(q)$ directed radially. These vectors are
perpendicular to each other, their scalar product hence is zero.
\par
    Thus, both summands in right hand side of formula \thetag{3.5}
vanish as $t\to 0$, hence the relationship \thetag{3.3} holds.
Similar to \thetag{3.2}, this relationship is fulfilled due to
initial data \thetag{2.13} regardless to the choice of force field
$\bold F$ of the dynamical system \thetag{1.1}. Therefore we consider
analogous relationship for second order derivatives
$$
\ddot\varphi_i\,\hbox{\vrule height 8pt depth 8pt width 0.5pt}_{\,t=0}
=\lim_{t\to 0}\frac{\partial^{\kern 0.5pt 2}\varphi_i(u^1,\ldots,u^{n-1},t)}
{\partial t^2}=0,\hskip -2em
\tag3.11
$$
which also should be fulfilled in the case of normal blow-up of
point $p_0$. Let's differentiate the equality \thetag{3.5} with
respect to $t$. This yields
$$
\ddot\varphi_i=(\nabla_t\bold F\,|\,\boldsymbol\tau_i)+2\,(\bold F\,|\,
\nabla_t\boldsymbol\tau_i)+(\bold v\,|\,\nabla_{tt}\boldsymbol\tau_i).
\hskip -2em
\tag3.12
$$
Components of force vector $\bold F$ depend on double set of arguments:
on coordinates $x^1,\,\ldots,\,x^n$ of the point on trajectory and on
components $v^1,\,\ldots,\,v^n$ of velocity vector of this point.
This means that vector $\bold F$ depend on the point of tangent bundle
$TM$. Such vectors are not embraced by the ordinary concept of vector
field on a manifold. Therefore in paper \cite{6} the concept of
{\bf extended vector field} was introduced. The concept of {\bf extended
tensor field} is its natural generalization.
\definition{Definition 2.3} {\bf Extended} tensor field $\bold X$
of the type $(r,s)$ on the manifold $M$ is a tensor-valued function
that to each point $q=(p,\bold v)$ of tangent bundle $TM$ puts into
correspondence some tensor from tensor space $T^r_s(p,M)$ at the point
$p$ on $M$.
\enddefinition
    Smooth extended tensor fields constitute an algebra over the ring
of smooth functions on tangent bundle. It was called {\bf the extended
algebra of tensor fields} on $M$. In the case of Riemannian manifold
one can naturally define two covariant differentiations $\nabla$ and
$\tilde\nabla$ in extended algebra of tensor fields on it. First was
called {\bf spatial gradient}, second was called {\bf velocity gradient}.
Covariant derivative $\nabla_t\bold F$ of the force vector in formula
\thetag{3.12} can be expressed through corresponding gradients of
extended vector field $\bold F$. For the components of the vector
$\nabla_t\bold F$ in formula \thetag{3.12} we have the following
expression:
$$
\nabla_t F^k=\sum^n_{s=1}\nabla_sF^k\,v^s+\sum^n_{s=1}\tilde\nabla_sF^k
\,F^s.\hskip -2em
\tag3.13
$$
We shall not comment formula \thetag{3.13}, and we shall not describe
in details all things related with extended algebra of tensor fields
(see Chapters \uppercase\expandafter{\romannumeral 2}, \uppercase
\expandafter{\romannumeral 3}, and \uppercase\expandafter{\romannumeral
4} in thesis \cite{17}). Technique of using extended tensor fields
is assumed to be known to reader.\par
     Second covariant derivative $\nabla_{tt}\boldsymbol\tau_i$ in
formula \thetag{3.12} is expressed through $\boldsymbol\tau_i$ and
$\nabla_t\boldsymbol\tau_i$. The matter is that components of any
vector of variation of trajectories $\boldsymbol\tau$ in case of
Newtonian dynamical systems satisfy the system of linear ordinary
differential equations of the second order:
$$
\aligned
\nabla_{tt}\tau^k&=-\sum^n_{m=1}\sum^n_{s=1}\sum^n_{r=1}R^k_{msr}\,
\tau^s\,v^r\,v^m+\\
&+\sum^n_{s=1}\nabla_t\tau^s\,\tilde\nabla_sF^k+\sum^n_{s=1}\tau^s\,
\nabla_sF^k.
\endaligned\hskip -2em
\tag3.14
$$
Taking into account \thetag{3.13} and \thetag{3.14}, we can bring
formula \thetag{3.12} to the form
$$
\gathered
\ddot\varphi_i=\sum^n_{r=1}\left(2\,F_r+\shave{\sum^n_{s=1}} v^s\,
\tilde\nabla_rF_s\right)\nabla_t\tau^r_i\,+\hskip -2em\\
+\,\sum^n_{r=1}\left(\,\shave{\sum^n_{s=1}}v^s\left(\nabla_sF_r+
\nabla_rF_s\right)+\shave{\sum^n_{s=1}} F^s\,\tilde\nabla_sF_r
\right)\tau^r_i.\hskip -2em
\endgathered
\tag3.15
$$
Here and everywhere below, aside with contravariant components of
vectors, we use their covariant components obtained by lowering
index by means of metric: 
$$
\xalignat 2
&v_i=\sum^n_{j=1}g_{ij}\,v^j,&&F_i=\sum^n_{j=1}g_{ij}\,F^j.
\endxalignat
$$
Quantities $F_k$, $\tilde\nabla_kF_s$, $\nabla_sF_k$ in formula
\thetag{3.15} are the components of smooth extended tensor fields.
They all have finite limits as $t\to 0$. Limits are determined by
substituting local coordinates of the point $p_0$ and components
of the vector $\bold v(0)=\nu_0\cdot\bold n(q)$ for their arguments.
Therefore in order to calculate limit of the derivative $\ddot
\varphi_i$ it is sufficient to use the relationship \thetag{3.9}
and remember that $\boldsymbol\tau_i\to 0$ as $t\to 0$ (the latter
is due to the expansions \thetag{2.16}):
$$
\lim_{t\to 0}\ddot\varphi_i=\sum^n_{r=1}\nu_0\,\left(2\,F_r
+\shave{\sum^n_{s=1}}v^s\,\tilde\nabla_rF_s\right)\,K^r_i.
\hskip -2em
\tag3.16
$$
Substituting \thetag{3.16} into \thetag{3.11}, we obtain the
following relationship:
$$
\sum^n_{r=1}\left(2\,F_r+\shave{\sum^n_{s=1}}v^s\,\tilde\nabla_rF_s
\right)K^r_i=0.\hskip -2em
\tag3.17
$$
Here $K^r_i$ are the components of the vector $\bold K_i$ from
\thetag{3.8}. Note that left hand side of \thetag{3.17} is linear
with respect to components of the vector $\bold K_i$, while vectors
$\bold K_1,\,\ldots,\,\bold K_{n-1}$ form a base in the hyperplane
perpendicular to the vector $\bold n(q)$. Vector $\bold n(q)$, in
turn, is collinear to the velocity vector $\bold v(0)=\nu_0\cdot
\bold n(q)$. Therefore if we introduce the operator $P$ of orthogonal
projection to the hyperplane perpendicular to velocity vector $\bold v$
and if we denote by $P^r_i$ its components, we can replace \thetag{3.17}
by an equivalent relationship 
$$
\sum^n_{r=1}\left(2\,F_r+\shave{\sum^n_{s=1}}v^s\,\tilde\nabla_rF_s
\right)P^r_i=0.\hskip -2em
\tag3.18
$$
Orthogonal projectors $\bold P$ form an extended tensor field of the
type $(1,1)$. Components of this field can be written in explicit form:
$$
P^r_i=\delta^r_i-N^r\,N_i.\hskip -2em
\tag3.19
$$
Here $\delta^r_i$ is Kronecker delta symbol, while $N^r$ and $N_i$
are contravariant and covariant components of extended vector field
$\bold N$ formed by unitary vectors collinear to the vector of
velocity $\bold v$:
$$
\xalignat 2
&v=|\bold v|,
&&\bold N=\frac{\bold v}{v}.\hskip -2em
\tag3.20
\endxalignat
$$
What is the meaning of the derived relationships \thetag{3.18}\,?
The matter is that the relationships \thetag{3.11}, in contrast to
\thetag{3.2} and \thetag{3.3}, cannot be fulfilled only due to
initial conditions \thetag{2.13}. They are equivalent to the
relationships \thetag{3.18} that should be fulfilled at the point
$p_0$ for all vectors $\bold v$ such that $|\bold v|=\nu_0$. If
dynamical system \thetag{1.1} belongs to the class of systems
admitting normal blow-up of points, as we assumed above in the
beginning of this section, then the relationships \thetag{3.18}
for its force field $\bold F$ are fulfilled at all points of
tangent bundle $TM$, where $|\bold v|\neq 0$. In this case
they are partial differential equations with respect to the
components of force vector $\bold F$.\par
     Further we continue to study the relationships like \thetag{3.11}.
Next in the series of relationships \thetag{3.2}, \thetag{3.3}, and
\thetag{3.11} is the vanishing condition for third derivatives of
the functions of deviation $\varphi_1,\,\ldots,\,\varphi_{n-1}$:
$$
\dddot\varphi_i\,\hbox{\vrule height 8pt depth 8pt width 0.5pt}_{\,t=0}
=\lim_{t\to 0}\frac{\partial^{\kern 0.5pt 3}\varphi_i(u^1,\ldots,u^{n-1},t)}
{\partial t^3}=0.\hskip -2em
\tag3.21
$$
In order to calculate third derivative $\dddot\varphi_i$ we differentiate
the equality \thetag{3.15} with respect to $t$. The equality \thetag{3.15}
has the following structure:
$$
\ddot\varphi_i=\sum^n_{r=1}\alpha_r\,\nabla_t\tau^r_i
+\sum^n_{r=1}\beta_r\,\tau^r_i.\hskip -2em
\tag3.22
$$
Here $\alpha_r$ and $\beta_r$ are components of extended covector fields.
Therefore 
$$
\gather
\dddot\varphi_i=\sum^n_{r=1}\alpha_r\,\nabla_{tt}\tau^r_i
+\sum^n_{r=1}\left(\nabla_t\alpha_r+\beta_r\right)\,\nabla_t\tau^r_i
+\sum^n_{r=1}\nabla_t\beta_r\,\tau^r_i=\\
=\sum^n_{r=1}\alpha_r\,\nabla_{tt}\tau^r_i+\sum^n_{r=1}
\left(\,\shave{\sum^n_{s=1}\nabla_s\alpha_r\,v^s+\sum^n_{s=1}
\tilde\nabla_s\alpha_r\,F^s}\right)\nabla_t\tau^r_i\,+\\
+\,\sum^n_{r=1}\beta_r\,\nabla_t\tau^r_i+\sum^n_{r=1}
\left(\,\shave{\sum^n_{s=1}\nabla_s\beta_r\,v^s+\sum^n_{s=1}
\tilde\nabla_s\beta_r\,F^s}\right)\tau^r_i.
\endgather
$$
In order to calculate $\nabla_t\alpha_r$ and $\nabla_t\beta_r$
above we used formulas similar to \thetag{3.13}. Further we
take into account that $\tau^r_i(0)=0$ (this follows from
\thetag{2.16}). Then 
$$
\dddot\varphi_i(0)=\sum^n_{r=1}\alpha_r\,\nabla_{tt}\tau^r_i(0)
+\sum^n_{r=1}\left(\beta_r+\shave{\sum^n_{s=1}\nabla_s\alpha_r\,v^s
+\sum^n_{s=1}\tilde\nabla_s\alpha_r\,F^s}\right)\nabla_t\tau^r_i(0).
$$
Second covariant derivative $\nabla_{tt}\tau^r_i$ can be determined
from the equation \thetag{3.14}. In the limit as $t\to 0$ this equation
yields
$$
\nabla_{tt}\tau^r_i(0)=\sum^n_{s=1}\tilde\nabla_sF^r
\,\nabla_t\tau^s_i(0).\hskip -2em
\tag3.23
$$
Let's substitute \thetag{3.23} into the above expression for
$\dddot\varphi_i(0)$. As a result we get
$$
\dddot\varphi_i(0)=\sum^n_{r=1}\left(\beta_r+\shave{\sum^n_{s=1}
\tilde\nabla_rF^s\,\alpha_s+\sum^n_{s=1}\nabla_s\alpha_r\,v^s
+\sum^n_{s=1}\tilde\nabla_s\alpha_r\,F^s}\right)\nabla_t\tau^r_i(0).
$$
The value of $\nabla_t\tau^r_i$ for $t=0$ is determined from \thetag{3.9}.
Therefore the condition of vanishing of third derivatives \thetag{3.21}
leads to the following relationship:
$$
\sum^n_{r=1}\left(\beta_r+\shave{\sum^n_{s=1}
\tilde\nabla_rF^s\,\alpha_s+\sum^n_{s=1}\nabla_s\alpha_r\,v^s
+\sum^n_{s=1}\tilde\nabla_s\alpha_r\,F^s}\right)K^r_i=0.
$$
Components of vectors $\bold K_1,\,\ldots,\,\bold K_{n-1}$ in the
relationship just obtained can be replaced by components of orthogonal
projector $\bold P$. Arguments for doing this are the same as in
replacing the relationship \thetag{3.17} by \thetag{3.18}:
$$
\sum^n_{r=1}\left(\beta_r+\shave{\sum^n_{s=1}\tilde\nabla_rF^s
\,\alpha_s+\sum^n_{s=1}\nabla_s\alpha_r\,v^s+\sum^n_{s=1}
\tilde\nabla_s\alpha_r\,F^s}\right)P^r_i=0.\hskip -2em
\tag3.24
$$
Now we are to substitute explicit expressions for $\alpha_r$ and
$\beta_r$ into the relationship \thetag{3.24}. They should be
taken in comparing formulas \thetag{3.15} and \thetag{3.22} for
$\ddot\varphi$:
$$
\aligned
&\alpha_r=2\,F_r+\sum^n_{s=1} v^s\,\tilde\nabla_rF_s,\\
&\beta_r=\sum^n_{s=1}v^s\left(\nabla_sF_r+
\nabla_rF_s\right)+\sum^n_{s=1}F^s\,\tilde\nabla_sF_r.
\endaligned\hskip-2em
\tag3.25
$$
But before doing this substitution, note that previously obtained
equations \thetag{3.18} for the components of force vector $\bold F$
can be written as 
$$
\sum^n_{r=1}\alpha_r\,P^r_i=0.\hskip-2em
\tag3.26
$$
Let's apply the differentiations $\nabla$ and $\tilde\nabla$ to
\thetag{3.26} and let's contract the resulting equalities with
components of vectors $\bold v$ and $\bold F$ respectively.
This yields
$$
\aligned
&\sum^n_{r=1}\sum^n_{s=1}v^s\,\nabla_s\alpha_r\,P^r_i=0,\\
&\sum^n_{r=1}\sum^n_{s=1}F^s\,\tilde\nabla_s\alpha_r\,P^r_i=
\sum^n_{r=1}\sum^n_{s=1}\frac{\alpha_s\,N^s\,P^r_i\,F_r}
{|\bold v|}.\endaligned\hskip -2em
\tag3.27
$$
In deriving \thetag{3.27} we took into account \thetag{3.26} and
we used the relationships
$$
\xalignat 2
&\quad\nabla_sP^r_i=0,&&\tilde\nabla_sP^r_i=-\frac{1}{|\bold v|}
\left(\shave{N_i\,P^r_s+\sum^n_{j=1}g_{ij}\,P^j_s\,N^r}\right).
\hskip -2em
\tag3.28
\endxalignat
$$
The relationships \thetag{3.28} can be proved by direct calculations
on the base of formulas \thetag{3.19} and \thetag{3.20} (see \S\,5
in Chapter~\uppercase\expandafter{\romannumeral 5} of thesis \cite{17}).
\par
    Let's use the relationships \thetag{3.27} in order to simplify the
equations \thetag{3.24}. Due to the first of these relationships the
third summand in \thetag{3.24} vanishes. Second relationship \thetag{3.27}
enables us to transform fourth summand in \thetag{3.24}. As a result of both
these transformations we obtain
$$
\sum^n_{r=1}\left(\beta_r+\shave{\sum^n_{s=1}\tilde\nabla_rF^s
\,\alpha_s+\sum^n_{s=1}\frac{\alpha_s\,N^s\,F_r}{|\bold v|}}
\right)P^r_i=0.\hskip -2em
\tag3.29
$$
Now let's substitute $\alpha_r$ and $\beta_r$ taken from \thetag{3.25}
into the equations \thetag{3.29}. Then
$$
\gather
\sum^n_{r=1}\sum^n_{s=1}v^s(\nabla_sF_r+\nabla_rF_s)\,P^r_i+
\sum^n_{r=1}\sum^n_{s=1}F^s\tilde\nabla_sF_r\,P^r_i\,+\hskip -2em\\
+\,\sum^n_{r=1}\sum^n_{s=1}2\,F_s\,\tilde\nabla_rF^s\,P^r_i
+\sum^n_{r=1}\sum^n_{s=1}\sum^n_{q=1}v^q\,\tilde\nabla_sF_q\,
\tilde\nabla_rF^s\,P^r_i\,+\hskip -2em
\tag3.30\\
+\,\sum^n_{r=1}\sum^n_{s=1}\frac{2\,F_r\,N^s\,F_s}{|\bold v|}\,P^r_i
+\sum^n_{r=1}\sum^n_{s=1}\sum^n_{q=1}F_r\,N^s\,N^q\,
\tilde\nabla_s\,F_q\,P^r_i=0.\hskip -2em
\endgather
$$
Further transformation of the obtained equations we begin with the
fourth summand in \thetag{3.30}. Due to \thetag{3.19} we have
$\delta^r_i=P^r_i+N^r\,N_i$. Therefore
$$
\gather
\sum^n_{r=1}\sum^n_{s=1}\sum^n_{q=1}v^q\,\tilde\nabla_sF_q\,
\tilde\nabla_rF^s\,P^r_i=\sum^n_{r=1}\sum^n_{s=1}\sum^n_{q=1}\sum^n_{j=1}
v^q\,\tilde\nabla_sF_q\,P^s_j\,\tilde\nabla_rF^j\,P^r_i\,+\\
+\sum^n_{r=1}\sum^n_{s=1}\sum^n_{q=1}\sum^n_{j=1}v^q\,\tilde\nabla_sF_q
\,N^s\,N_j\,\tilde\nabla_rF^j\,P^r_i=-\sum^n_{r=1}\sum^n_{s=1}
\sum^n_{j=1}2\,F_s\,P^s_j\,\tilde\nabla_rF^j\,P^r_i\,+\\
+\sum^n_{r=1}\sum^n_{s=1}\sum^n_{q=1}\sum^n_{j=1}N^q\,\tilde\nabla_sF_q
\,N^s\,v^j\,\tilde\nabla_rF_j\,P^r_i=-\sum^n_{r=1}\sum^n_{s=1}
2\,F_s\,\tilde\nabla_rF^s\,P^r_i\,+\\
+\sum^n_{r=1}\sum^n_{s=1}\sum^n_{j=1}2\,F_j\,N^j\,N_s\,
\tilde\nabla_rF^s\,P^r_i-\sum^n_{r=1}\sum^n_{s=1}\sum^n_{q=1}
2\,N^q\,N^s\,\tilde\nabla_sF_q\,F_r\,P^r_i.
\endgather
$$
Apart from the relationship $\delta^r_i=P^r_i+N^r\,N_i$ following
from \thetag{3.19}, here we used the relationship \thetag{3.26}
written as 
$$
\sum^n_{r=1}\sum^n_{s=1} |\bold v|\,N^s\,\tilde\nabla_rF_s\,P^r_i=-
\sum^n_{r=1}2\,F_r\,P^r_i.\hskip-2em
$$
Let's apply this relationship once more for to transform second
summand in the above expression. As a result we have
$$
\gather
\sum^n_{r=1}\sum^n_{s=1}\sum^n_{q=1}v^q\,\tilde\nabla_sF_q\,
\tilde\nabla_rF^s\,P^r_i=-\sum^n_{r=1}\sum^n_{s=1}
2\,F_s\,\tilde\nabla_rF^s\,P^r_i\,-\\
-\sum^n_{r=1}\sum^n_{s=1}\frac{4\,F_s\,N^s\,F_r}{|\bold v|}
\,P^r_i-\sum^n_{r=1}\sum^n_{s=1}\sum^n_{q=1}
2\,N^q\,N^s\,\tilde\nabla_sF_q\,F_r\,P^r_i.
\endgather
$$
In substituting this expression for the fourth summand into
\thetag{3.30} we find that third summand cancels, while
fifth and sixth summands change their signs. The equations
\thetag{3.30} in whole now look like 
$$
\gathered
\sum^n_{r=1}\sum^n_{s=1}v^s(\nabla_sF_r+\nabla_rF_s)\,P^r_i+
\sum^n_{r=1}\sum^n_{s=1}F^s\tilde\nabla_sF_r\,P^r_i\,-\\
-\,\sum^n_{r=1}\sum^n_{s=1}\frac{2\,F_r\,N^s\,F_s}{|\bold v|}\,P^r_i
-\sum^n_{r=1}\sum^n_{s=1}\sum^n_{q=1}F_r\,N^s\,N^q\,
\tilde\nabla_s\,F_q\,P^r_i=0.
\endgathered\hskip -2em
\tag3.31
$$
And finally, let's do some slight (purely cosmetic) transformations
in the equations \thetag{3.18} and \thetag{3.31}. Then write them
combining into a system:
$$
\cases
\dsize\sum^n_{i=1}\left(v^{-1}\,F_i+\shave{\sum^n_{j=1}}
\tilde\nabla_i\left(N^j\,F_j\right)\right)P^i_k=0,\\
\vspace{2ex}
\aligned
 &\sum^n_{i=1}\sum^n_{j=1}\left(\nabla_iF_j+\nabla_jF_i-2\,v^{-2}
 \,F_i\,F_j\right)N^j\,P^i_k\,+\\
 &+\sum^n_{i=1}\sum^n_{j=1}\left(\frac{F^j\,\tilde\nabla_jF_i}{v}
 -\sum^n_{r=1}\frac{N^r\,N^j\,\tilde\nabla_jF_r}{v}\,F_i\right)P^i_k=0.
 \endaligned
\endcases\hskip -2em
\tag3.32
$$
Let's state the result following from the above calculations
in form of a theorem.
\proclaim{Theorem 3.1} For Newtonian dynamical system \thetag{1.1}
in Riemannian manifold $M$ to admit normal blow-up of points its
force field should satisfy the equations \thetag{3.32} at all points
$q=(p,\bold v)$ of tangent bundle $TM$, where $\bold v\neq 0$.
\endproclaim
\head
4. Weak normality condition.
\endhead
    Above we considered vanishing conditions for the functions
of deviation $\varphi_i$ and for their derivatives $\dot\varphi_i$,
$\ddot\varphi_i$, and $\dddot\varphi_i$ at the initial instant of time
$t=0$. And we have found that first two conditions $\varphi_i(0)=0$
and $\dot\varphi_i(0)=0$ are fulfilled only due to initial data
\thetag{2.13} determining the blow-up of point. They make no
restriction for the choice of force field $\bold F$ of dynamical
system. Considering next two conditions $\ddot\varphi_i(0)=0$ and
$\dddot\varphi_i(0)=0$, we derived the restrictions for $\bold F$
in form of the equations \thetag{3.32} for the components of force
vector.\par
    Further we could step by step consider the vanishing conditions
for the derivatives of the functions of deviation of higher order
getting more and more equations for $\bold F$ in each step. However,
as we shall see soon, it is not necessary. The matter is that the
equations \thetag{3.32} are exactly the same as {\bf weak normality
equations}, which arose in considering Newtonian dynamical systems
admitting the normal shift of hypersurfaces. For the case $M=\Bbb R^2$
they were first derived in \cite{1}, then in \cite{3} they were
generalized for the case $M=\Bbb R^n$. In form of \thetag{3.32}
corresponding to the case of arbitrary Riemannian manifold these
equations were derived in \cite{6} (see also Chapter \uppercase
\expandafter{\romannumeral 5} in thesis \cite{17}). Weak normality
equations are equivalent to the following condition {\bf of weak
normality}.
\definition{Definition 4.1} Newtonian dynamical system \thetag{1.1}
in Riemannian manifold $M$ satisfies {\bf weak normality condition}
if for each its trajectory there exists some ordinary differential
equation
$$
\ddot\varphi=\Cal A(t)\,\dot\varphi+\Cal B(t)\,\varphi\hskip -2em
\tag4.1
$$
such that any function of deviation $\varphi=(\bold v\,|\,\boldsymbol
\tau)$ corresponding to the arbitrary choice of the vector of variation 
$\boldsymbol\tau$ on that trajectory is the solution of this equation.
\enddefinition
     Words ``any function of deviation'' and ``arbitrary choice of the
vector of variation'' in this definition should be commented. Suppose
that $p=p(t)$ is some trajectory of Newtonian dynamical system
\thetag{1.1}. Let's include it into some (arbitrary) one-parametric
family of trajectories $p=p(u,t)$, so that for $u=0$ we would have
$p(0,t)=p(t)$. Parameter $u$ can be introduced, for instance, by making
dependent on $u$ the initial data in Cauchy problem that fixes our
trajectory $p(t)$. In local coordinates the family of trajectories
$p=p(u,t)$ is given by the functions
$$
\align
&x^1=x^1(u,t),\hskip -2em\\
&.\ .\ .\ .\ .\ .\ .\ .\ .\ .\ 
\hskip -2em
\tag4.2\\
&x^n=x^n(u,t).\hskip -2em
\endalign
$$
Derivatives of the functions \thetag{4.2} with respect to parameter $u$
determine vector $\boldsymbol\tau(t)$ on the trajectory $p(t)$:
$$
\boldsymbol\tau=\sum^n_{k=1}\frac{\partial x^k}{\partial u}\,
\hbox{\vrule height 12pt depth 8pt width 0.5pt}_{\,u=0}\cdot
\frac{\partial}{\partial x^k}\hskip -2em
\tag4.3
$$
(compare with formula \thetag{2.6}). This vector is called {\bf the
vector of variation of trajectories}. This vector $\boldsymbol\tau(t)$
(constructed as described above) is implied in definition~4.1. Its
scalar product with the vector of velocity $\bold v$ is {\bf a function
of deviation} corresponding to it: $\varphi=(\bold v\,|\,\boldsymbol
\tau)$.\par
     It is easy to show that components of any vector of variation
$\boldsymbol\tau(t)$ constructed as above satisfy the differential
equations \thetag{3.14} (see paper \cite{6} or Chapter \uppercase
\expandafter{\romannumeral 5} in thesis \cite{17}). And conversely,
any vector $\boldsymbol\tau(t)$ with components satisfying the
equations \thetag{3.14} can be obtained by formula \thetag{4.3} in
the above construction. Therefore words ``arbitrary choice of the
vector of variation'' in definition~4.1 can be understood as the
choice of an arbitrary solution of the system of linear homogeneous
ordinary differential equations \thetag{3.14}. Function of deviation
$$
\varphi(t)=(\bold v\,|\,\boldsymbol\tau)=\sum^n_{k=1}v_k\,\tau^k
$$
corresponding to this choice of $\boldsymbol\tau$ in the case of
general position satisfies linear homogeneous ordinary differential
equation of the order $2\,n$ (see theorem~6.1 in Chapter \uppercase
\expandafter{\romannumeral 5} of thesis \cite{17}). In special cases
(for special choice of force field $\bold F$) the order of this
equation can be lower. Definition~4.1 separates the case, when order
of the equation for $\varphi$ is $2$. For this case in \cite{6}
the following proposition was proved (see also theorem~6.2 in Chapter
\uppercase\expandafter{\romannumeral 5} of thesis \cite{17}).
\proclaim{Theorem 4.1} Newtonian dynamical system in Riemannian manifold
of the dimension $n\geqslant 2$ satisfies weak normality condition if and
only if its force field satisfies the system of differential equations
\thetag{3.32} at all points $q=(p,\bold v)$ of tangent bundle $TM$,
where $\bold v\neq 0$.
\endproclaim
    Note that vectors of variation $\boldsymbol\tau_1,\,\ldots,\,
\boldsymbol\tau_{n-1}$ arising in blow-up of points are naturally
embedded into the above construction \thetag{4.3}. Therefore if
force field $\bold F$ of dynamical system \thetag{1.1} satisfies
the equations \thetag{3.32}, then corresponding functions of
deviation $\varphi_1,\,\ldots,\,\varphi_{n-1}$ satisfy the
differential equations of the form \thetag{4.1}. In this case the
conditions 
$$
\xalignat 2
&\varphi_i\,\hbox{\vrule height 8pt depth 8pt width 0.5pt}_{\,t=0}=0,
&&\dot\varphi_i\,\hbox{\vrule height 8pt depth 8pt width 0.5pt}_{\,t=0}
=0\hskip -2em
\tag4.4
\endxalignat
$$
from \thetag{3.2} and \thetag{3.3} provide identical vanishing of the
functions of deviation. The conditions \thetag{4.4} by themselves, as
we noted above, are provided only by initial data \thetag{2.13}.
Therefore we can strengthen the theorem~3.1 formulating it as follows.
\proclaim{Theorem 4.2} Newtonian dynamical system \thetag{1.1}
in Riemannian manifold $M$ admits normal blow-up of points if and only
if its force field satisfies the equations \thetag{3.32} at all points
$q=(p,\bold v)$ of tangent bundle $TM$, where $\bold v\neq 0$.
\endproclaim
    The restriction in dimension $n\geqslant 2$ from theorem~4.1 is
inessential. We do not formulate it explicitly, since in the dimension
$n=1$ the concept of normal blow-up of points has no meaning.\par
\head
5. Concluding remarks.
\endhead
    Theorem~4.2 is a central result of present paper. It reduces the
study of Newtonian dynamical systems admitting normal blow-up of points
to the analysis of the system of partial differential equations
\thetag{3.32} for their force fields. Moreover it provides geometric
interpretation for weak normality condition, reducing this rather abstract
condition to visually obvious condition that dynamical system is able
to perform normal blow-up of points.\par
    Note that in the dimension $n=2$ (this case was studied in details
in thesis \cite{18}) the condition that Newtonian dynamical system
admits the normal shift of hypersurfaces is also reduced to the system
of the equations \thetag{3.32}. Therefore we have the following
proposition.
\proclaim{Proposition 5.1} Class of Newtonian dynamical systems
admitting the normal shift of hypersurfaces and class of Newtonian
dynamical systems admitting normal blow-up of points for $n=2$
do coincide.
\endproclaim
    For the dimension $n\geqslant 3$ in the theory of Newtonian
dynamical systems admitting the normal shift of hypersurfaces, apart from
weak normality equations \thetag{3.32}, we have so called {\bf additional
normality equations}. They have the following form:
$$
\cases
\aligned
&\sum^n_{i=1}\sum^n_{j=1}P^i_\varepsilon\,P^j_\sigma\left(
\,\shave{\sum^n_{m=1}}N^m\,\frac{F_i\,\tilde\nabla_mF_j}{v}
-\nabla_iF_j\right)=\\
&\quad=\sum^n_{i=1}\sum^n_{j=1}P^i_\varepsilon\,P^j_\sigma\left(
\,\shave{\sum^n_{m=1}}N^m\,\frac{F_j\,\tilde\nabla_mF_i}{v}
-\nabla_jF_i\right),
\endaligned\\
\vspace{2ex}
\dsize\sum^n_{i=1}\sum^n_{j=1}P^j_\sigma\,\tilde\nabla_jF^i\,
P^\varepsilon_i=\sum^n_{i=1}\sum^n_{j=1}\sum^n_{m=1}
\frac{P^j_m\,\tilde\nabla_jF^i\,P^m_i}{n-1}\,P^\varepsilon_\sigma.
\endcases\hskip -2em
\tag5.1
$$
Analysis of complete system of normality equations combined by
\thetag{3.32} and \thetag{5.1} was undertook in \cite{16}. However,
in paper \cite{16} an error was committed. Therefore part of results
of \cite{16} are not valid. This error was corrected in Chapter
\uppercase\expandafter{\romannumeral 7} of thesis \cite{17} (see also
paper \cite{19}). As a result an explicit formula for general solution
of complete system of normality equations \thetag{3.32} and \thetag{5.1}
was derived.\par
    Currently the analysis of separate system of weak normality equations
\thetag{3.32} is urgent. In particular, would be worth to know whether
something like proposition~5.1 is valid in the dimension $n\geqslant 3$.
Theorem~4.2 reduces this problem to the study of the equations \thetag{3.32}
\head
6. Acknowledgments.
\endhead
     This work is supported by grant from Russian Fund for Basic Research
(project No\nolinebreak\.~00\nolinebreak-01-00068, coordinator 
Ya\.~T.~Sultanaev), and by grant from Academy of Sciences of the
Republic Bashkortostan (coordinator N.~M.~Asadullin). Author is grateful
to these organizations for financial support.


\Refs
\ref\no 1\by Boldin~A.~Yu\., Sharipov~R.~A.\book Dynamical systems
accepting the normal shift\publ Preprint No\.~0001-M of Bashkir State
University\publaddr Ufa\yr April, 1993
\endref
\ref\no 2\by Boldin~A.~Yu\., Sharipov~R.~A.\paper Dynamical systems
accepting the normal shift\jour Theoretical and Mathematical Physics (TMF)
\vol 97\issue 3\yr 1993\pages 386--395\moreref see also chao-dyn/9403003
in Electronic Archive at LANL\footnotemark
\endref
\footnotetext{Electronic Archive at Los Alamos national Laboratory of USA
(LANL). Archive is accessible through Internet 
{\bf http:/\negskp/xxx\.lanl\.gov}, it has mirror site 
{\bf http:/\negskp/xxx\.itep\.ru} at the Institute for Theoretical and
Experimental Physics (ITEP, Moscow).}\adjustfootnotemark{-1}
\ref\no 3\by Boldin~A.~Yu\., Sharipov~R.~A.\paper Multidimensional
dynamical systems accepting the normal shift\jour Theoretical and
Mathematical Physics (TMF)\vol 100\issue 2\yr 1994\pages 264--269
\moreref see also patt-sol/9404001 in Electronic Archive at LANL
\endref
\ref\no 4\by Boldin~A.~Yu\., Sharipov~R.~A.\paper Dynamical systems
accepting the normal shift\jour Reports of Russian Academy of Sciences
(Dokladi RAN)\vol 334\yr 1994\issue 2\pages 165--167
\endref
\ref\no 5\by Sharipov~R.~A.\paper Problem of metrizability for
the dynamical systems accepting the normal shift\jour Theoretical and
Mathematical Physics (TMF)\yr 1994\vol 101\issue 1\pages 85--93\moreref
see also solv-int/9404003 in Electronic Archive at LANL
\endref
\ref\no 6\by Boldin~A.~Yu\., Dmitrieva~V.~V., Safin~S.~S., Sharipov~R.~A.
\paper Dynamical systems accepting the normal shift on an arbitrary 
Riemannian manifold\jour Theoretical and Mathematical Physics (TMF)
\yr 1995\vol 105\issue 2\pages 256--266\moreref\inbook see also
``{Dynamical systems accepting the normal shift}'', Collection of papers
\publ Bashkir State University\publaddr Ufa\yr 1994\pages 4--19
\moreref see also hep-th/9405021 in Electronic Archive at LANL
\endref
\ref\no 7\by Boldin~A.~Yu\., Bronnikov~A.~A., Dmitrieva~V.~V.,
Sharipov~R.~A.\paper Complete normality conditions for the dynamical
systems on Riemannian manifolds\jour Theoretical and Mathematical
Physics (TMF)\yr 1995\vol 103\issue 2\pages 267--275\moreref\inbook
see also ``{Dynamical systems accepting the normal shift}'', Collection
of papers\publ Bashkir State University\publaddr Ufa\yr 1994
\pages 20--30\moreref see also astro-ph/9405049 in Electronic Archive
at LANL
\endref
\ref\no 8\by Boldin~A\.~Yu\.\paper On the self-similar solutions of 
normality equation in two-dimensional case\inbook ``{Dynamical systems
accepting the normal shift}'', Collection of papers\publ Bashkir State
University\publaddr Ufa\yr 1994\pages 31--39\moreref see also
patt-sol/9407002 in Electronic Archive at LANL
\endref
\ref\no 9\by Sharipov~R.~A.\paper Metrizability by means of conformally
equivalent metric for the dynamical systems\jour Theoretical and
Mathematical Physics (TMF)\yr 1995\vol 105\issue 2\pages 276--282
\moreref\inbook see also ``{Integrability in dynamical systems}''\publ
Institute of Mathematics, Bashkir Scientific Center of Ural branch of
Russian Academy of Sciences (BNC UrO RAN)\publaddr Ufa\yr 1994
\pages 80--90
\endref
\ref\no 10\by Sharipov~R\.~A\.\paper Dynamical systems accepting normal
shift in Finslerian geometry,\yr November, 1993\finalinfo 
unpublished\footnotemark
\endref
\footnotetext{Papers \cite{1--16} are arranged here in the order they
were written. However, the order of publication not always coincides with
the order of writing.}
\ref\no 11\by Sharipov~R\.~A\.\paper Normality conditions and affine
variations of connection on Riemannian manifolds,\yr December, 1993
\finalinfo unpublished
\endref
\ref\no 12\by Sharipov~R.~A.\paper Dynamical system accepting the normal
shift (report at the conference)\jour see in Progress in Mathematical
Sciences (Uspehi Mat\. Nauk)\vol 49\yr 1994\issue 4\page 105
\endref
\ref\no 13\by Sharipov~R.~A.\paper Higher dynamical systems accepting 
the normal shift\inbook ``{Dynamical systems accepting the normal 
shift}'', Collection of papers\publ Bashkir State University\publaddr 
Ufa\yr 1994\linebreak\pages 41--65
\endref
\ref\no 14\by Dmitrieva~V.~V.\paper On the equivalence of two forms
of normality equations in $\Bbb R^n$\inbook ``{Integrability in dynamical
systems}''\publ Institute of Mathematics, Bashkir Scientific Center of
Ural branch of Russian Academy of Sciences (BNC UrO RAN)\publaddr
Ufa\yr 1994\pages 5--16
\endref
\ref\no 15\by Bronnikov~A.~A., Sharipov~R.~A.\paper Axially
symmetric dynamical systems accep\-ting the normal shift in $\Bbb R^n$
\inbook ``{Integrability in dynamical systems}''\publ Institute of
Mathematics, Bashkir Scientific Center of Ural branch of Russian Academy
of Sciences (BNC UrO RAN)\publaddr Ufa\yr 1994\linebreak\pages 62--69
\endref
\ref\no 16\by Boldin~A.~Yu\., Sharipov~R.~A.\paper On the solution
of normality equations in the dimension $n\geqslant 3$\jour Algebra and
Analysis (Algebra i Analiz)\vol 10\yr 1998\issue 4\pages 37--62\moreref
see also solv-int/9610006 in Electronic Archive at LANL
\endref
\ref\no 17\by Sharipov~R.~A.\book Dynamical systems admitting the normal
shift\publ Thesis for the degree of Doctor of Sciences in Russia\yr 1999
\moreref English version of thesis is submitted to Electronic Archive at 
LANL, see archive file math.DG/0002202 in the section of Differential 
Geometry\footnotemark
\endref
\footnotetext{For the convenience of reader we give direct reference
to archive file. This is the following URL address:
{\bf http:/\negskp/xxx\.lanl\.gov/eprint/math\.DG/0002202}\,.}
\ref\no 18\by Boldin~A.~Yu\.\book Two-dimensional dynamical systems
admitting the normal shift\publ Thesis for the degree of Candidate of
Sciences in Russia\yr 2000
\endref
\ref\no 19\by Sharipov~R.~A.\paper Newtonian normal shift in
multidimensional Riemannian geometry\jour Paper math.DG/0006125
in Electronic Archive at LANL\yr 2000
\endref
\endRefs
\enddocument
\end